\documentclass[12pt]{amsart}
\usepackage{a4wide,graphicx}

\def\Z{\mathbf Z}

\theoremstyle{plain}
\newtheorem{theorem}{Theorem}
\newtheorem{lemma}[theorem]{Lemma}

\newcommand{\burau}{\operatorname{Burau}}

\begin{document}

\title
    [The kernel of $\operatorname{Burau}(4)\otimes \Z_p$ is all pseudo-Anosov]
    {The kernel of {\boldmath $\operatorname{Burau}(4)\otimes \Z_p$}
is all pseudo-Anosov}
\author{Sang Jin Lee \and Won Taek Song}

\thanks{The research of the first author was supported
by the faculty research fund
of Konkuk University in 2003}
\address{Department of Mathematics, Konkuk University,
1 Hwayang-dong Gwangjin-gu, Seoul 143-701, Korea}
\email{sangjin@konkuk.ac.kr}

\address{School of Mathematics,
Korea Institute for Advanced Study,
207-43 Cheongnyangni 2-dong, Dongdaemun-gu,
Seoul 130-722, Korea}
\email{cape@kias.re.kr}

\date{February 11, 2004.}

\maketitle

\section{Introduction}
Given two pseudo-Anosov homeomorphisms with distinct invariant
measured foliations, some powers of their isotopy classes
generate a rank two free subgroup of the mapping class group
of the surface~\cite{MR87c:57009}.
This construction gives an example of all pseudo-Anosov subgroup
of the mapping class group.
Whittlesey~\cite{MR2001j:57022} gives a positive answer
to the natural question of the existence of all pseudo-Anosov
\emph{normal} subgroups by showing that the Brunnian mapping classes on
a sphere with at least five punctures are neither periodic nor
reducible.

In this note we show that the kernel of $\burau(4) \otimes \Z_p$,
the reduced Burau representation with coefficients in $\Z_p$ of
the 4-braid group $B_4$, consists only of pseudo-Anosov braids.
By Cooper and Long~\cite{MR97m:20050,MR99k:20077},
$\burau(4)\otimes\Z_p$ for $p=2,3$ is not faithful.
It is straightforward to check that there exist non-Brunnian braids
in the kernels, hence giving new examples of all pseudo-Anosov normal subgroups
of $B_4$ that are not contained in the example of Whittlesey.

For the proof, assume that we are given a
nontrivial 4-braid that is not pseudo-Anosov.
If it is periodic, it is conjugate to a rigid rotation~\cite{Brouwer1919},
whose Burau action is clearly non-trivial.
If it is reducible, then in many ways
it is similar to a 3-braid so that
its Burau action is fairly predictable, for which case
an automaton that records the polynomial degrees
suffices to prove faithfulness.
Our argument is similar to that of the ping-pong lemma.
We construct an automaton whose states are disjoint subsets
of $\Z_p[t,t^{-1}]^3$ and whose arrows are braid actions
that map the subsets into the subsets.

For braids with more than four strands,
this approach is immediately confronted by
various obstacles.
Since $\burau(4)\otimes \Z_2$ is not faithful,
the kernel of  $\burau(5)\otimes \Z_2$ contains
reducible braids.
Taking other representations or taking intersection with
other subgroups to get rid of such reducible braids
then makes the proof more difficult.

We remark that the present result  is a byproduct of working
on the faithfulness question
of $\burau(4)$~\cite{MR92b:20041,MR93k:57019,MR94c:20071,MR2001j:20055}.

\section{$\ker \burau(4)\otimes \Z_p$
does not contain periodic or reducible braids.}
The $n$-braid group $B_n$ consists of the mapping classes on
the $n$-punctured disk. The center of $B_n$ is the infinite
cyclic group generated by
the Dehn twist along the boundary.
A braid is called \emph{periodic} if
some of its power is contained in the center.
A braid  is called \emph{reducible} if  it is represented by
It  is called \emph{reducible} if  it is represented by
a disk homeomorphism that fixes a collection of disjoint essential curves.
If a braid is neither periodic nor reducible, then by  Nielsen-Thurston
classification of surface homeomorphisms~\cite{MR89k:57023,MR82m:57003} it is represented by a pseudo-Anosov
homeomorphism. Such a braid is called \emph{pseudo-Anosov}.
A subgroup of $B_n$ is called \emph{all pseudo-Anosov} if
its non-trivial elements are all pseudo-Anosov.

The $n$-braid group $B_n$ has the following presentation:
$$ B_n = \left\langle
\sigma_1,\ldots,\sigma_{n-1}\biggm|
\begin{array}{l}
\sigma_i\sigma_j=\sigma_j\sigma_i,\quad |i-j|\ge 2\\
\sigma_i\sigma_j\sigma_i=\sigma_j\sigma_i\sigma_j,\quad |i-j|=1
\end{array}
\right\rangle
$$
The reduced Burau representation $\rho_n = \burau(n)\colon
B_n \to GL_{n-1}(\Z[t,t^{-1}])$ is defined by the action on the first
homology of the cyclic cover of the punctured disk.
For the purpose of this note, it suffices to define
$\rho_4$ by the following three matrices.
$$ \rho_4(\sigma_1)= \left(\begin{array}{ccc}
-t & 0 & 0\\
 1 & 1 & 0\\
 0 & 0 & 1
\end{array}\right)\qquad
\rho_4(\sigma_2)= \left(\begin{array}{ccc}
1 &  t & 0\\
0 & -t & 0\\
0 &  1 & 1
\end{array}\right)\qquad
\rho_4(\sigma_3)= \left(\begin{array}{ccc}
1 & 0 & 0\\
0 & 1 & t\\
0 & 0 & -t\\
\end{array}\right)
$$
We use the convention that $B_4$ acts on $\Z[t,t^{-1}]^3$ from
the right. We denote by $\mathbf v*_\rho \beta$, or more simply by
$\mathbf v* \beta$,
the matrix multiplication
$\mathbf v\rho(\beta) $ for
a row vector $\mathbf v$, a representation $\rho$ and  a braid $\beta$.
For example,
$(f,g,h)*_{\rho_4} \sigma_1=(-tf+g,g,h)$ for $f,g,h\in\Z[t,t^{-1}]$.

The following theorem is the main result of this note. 
\begin{theorem}\label{thm:main}
The kernel of $\rho_4\otimes \Z_p\colon B_4 \to GL_3(\Z_p[t,t^{-1}])$ for
$p\ge 2$
does not contain  a non-trivial periodic or reducible  braid.
In particular if $\rho_4\otimes \Z_p$ is not faithful,
its kernel is an all pseudo-Anosov normal subgroup of $B_4$.
\end{theorem}

\begin{lemma}\label{thm:periodic} $\rho_n\otimes \Z_p$ is
faithful for periodic braids.
\end{lemma}
\begin{proof}
If $\beta\in B_n$ is a periodic $n$-braid, then
it is represented by a rigid rotation on the punctured disk~\cite{Brouwer1919}
so that it is conjugate to $(\sigma_{n-1}\cdots \sigma_2\sigma_1)^k$
or to $(\sigma_{n-1}\cdots\sigma_2\sigma_1\sigma_1)^k$ for some $k\in\Z$.
Since $\det((\rho_n\otimes\Z_p)(\beta))=(-t)^{e(\beta)}$,
where the exponent sum $e(\beta)$ is $k(n-1)$ or $kn$;
if $\beta$ is in the kernel of $\rho_n\otimes\Z_p$,
then $k=0$ and $\beta$ is trivial.
\end{proof}

Let $\Delta_3=\sigma_1\sigma_2\sigma_1 \in B_3$ and
$\Delta_4=\sigma_1\sigma_2\sigma_1\sigma_3\sigma_2\sigma_1\in B_4$
be the square roots of the generator of the center of $B_3$
and $B_4$, respectively.
For a Laurent polynomial $f(t)=\sum_m a_mt^m$,
define $\deg f=\max\{m:a_m\ne 0\}$.
By convention we define $\deg f=-\infty$ if $f=0$.

\begin{lemma}\label{thm:3braid}
$\rho_3\otimes \Z_p$ is faithful.
\end{lemma}
\begin{proof}
Let $\rho =\rho_3\otimes \Z_p$ be the reduced Burau representation of $B_3$
with coefficients in $\Z_p$.
It is given by the following matrices.
$$ \rho(\sigma_1)= \left(\begin{array}{cc}
-t & 0 \\
 1 & 1
\end{array}\right) \qquad
\rho(\sigma_2)= \left(\begin{array}{cc}
1 &  t \\
0 & -t
\end{array}\right)
$$

Suppose that $\rho(\beta)$ is trivial for some non-trivial 3-braid $\beta$.
By Lemma~\ref{thm:periodic}, it is either reducible or pseudo-Anosov.
If $\beta$ is
reducible, then it is conjugate to
$\Delta_3^{2m}\sigma_1^k(\sigma_2\sigma_1^2\sigma_2)^l$ for
some integers $k$, $l$ and $m$,
which is an arbitrary 3-braid with an invariant curve standardly embedded
in the disk enclosing the first two punctures
as in Figure~\ref{fig:red}~(b).
By the relation
$\sigma_2\sigma_1^2\sigma_2=\Delta_3^2\sigma_1^{-2}$,
the 3-braid $\beta$ is conjugate to $\Delta_3^{2l+2m}\sigma_1^{k-2l}$.
Since $\rho(\beta)$ is trivial,
$$\rho(\Delta_3^{2l+2m}\sigma_1^{k-2l})
=t^{3(l+m)}\left(
\begin{array}{cc}
(-t)^{k-2l} & 0\\
* & 1
\end{array}\right)
$$
must be the identity matrix. So $l+m=0$ and $k-2l=0$ hence
$\beta$ is trivial, which contradicts the assumption.

If $\beta$ is pseudo-Anosov, it is conjugate to
$P(\sigma_1^{-1},\sigma_2)\Delta_3^{2k}$
where $P$ is a positive word on two letters~\cite{MR50:8496,MR1915500}.
By taking inverse or conjugation by $\Delta_3$ if necessary,
we can assume that $P(\sigma_1^{-1},\sigma_2)$ starts with
$\sigma_2$.
In other words $\beta$ or $\beta^{-1}$ is conjugate to
$$\alpha=\sigma_2Q(\sigma_1^{-1},\sigma_2)\Delta_3^{2k}
$$
for some positive word $Q$.
The $\rho$-actions of $\sigma_1^{-1}$,
$\sigma_2$ and $\Delta_3^2$ on $\Z_p[t,t^{-1}]^2$ are given as follows:
for $\mathbf v=(f,g) \in \Z_p[t,t^{-1}]^2$,
$$\textstyle
\mathbf v* \sigma_1^{-1} =(-{t^{-1}}(f-g),g),\quad
\mathbf v*\sigma_2=(f,t(f-g))\quad\mbox{and}\quad
\mathbf v* \Delta_3^2 =(t^3f,t^3g).
$$
Consider the subset $V_0=\{(f,g)\in\Z[t,t^{-1}]^2\mid \deg f<\deg
g\}$. It is easy to check that $V_0$ is invariant
under the action of $\sigma_1^{-1}$, $\sigma_2$ and $\Delta_3^{2}$.
Let $\mathbf v_0=(1,0)$. Then $\mathbf v_0*\sigma_2=(1,t)\in V_0$
so that
$\mathbf v_0*\alpha=(1,t)*Q(\sigma_1^{-1},\sigma_2)\Delta_3^{2k}
\in V_0$. Since
$\mathbf v_0\not\in V_0$, we have  $\mathbf v_0*\alpha\ne\mathbf v_0$,
which contradicts the assumption that $\beta$ is in the kernel of~$\rho$.
\end{proof}

\begin{proof}[proof of Theorem~\ref{thm:main}]
Let $\rho =\rho_4\otimes \Z_p$ be the reduced Burau representation of $B_4$
with coefficients in $\Z_p$.
Assume $\rho(\beta)$ is trivial for
some non-trivial 4-braid $\beta \in B_4$.
The braid $\beta$ is either reducible or pseudo-Anosov
by Lemma~\ref{thm:periodic}.
We need to show that $\beta$ is not reducible.

Suppose that $\beta$ is reducible.
By taking some power of $\beta$ if necessary,
we may assume that $\beta$ is represented by a
homeomorphism that fixes an essential simple closed curve $C$.
By applying a conjugation by a braid that
sends $C$ to one of the curves in Figure~\ref{fig:red}, we
assume that $C$ is one of the two standardly embedded curves
and the homeomorphism representing $\beta$ fixes $C$.

\begin{figure}
$$
{\includegraphics{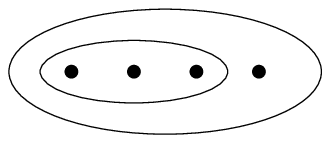}\atop\mbox{(a)}} \qquad
{\includegraphics{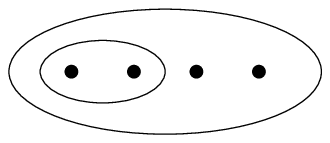}\atop\mbox{(b)}}
$$
\caption{Up to homeomorphisms on a 4-punctured disk,
there are only two essential curves.}\label{fig:red}
\end{figure}

Let $C$ be the curve enclosing the first three punctures as
Figure~\ref{fig:red}~(a).
Then $\beta$ can be written as
$\beta=\Delta_4^{2l}(\sigma_3\sigma_2\sigma_1^2\sigma_2\sigma_3)^k
W(\sigma_1,\sigma_2)$
for some integers $k$, $l$ and a word $W$ on two letters.
By the relation $ (\sigma_3\sigma_2\sigma_1^2\sigma_2\sigma_3) =
\Delta_4^2  (\sigma_1\sigma_2\sigma_1)^{-2}$,
we rewrite $\beta$ as
$$
\beta=\Delta_4^{2k+2l}(\sigma_1\sigma_2\sigma_1)^{-2k}W(\sigma_1,\sigma_2).
$$
Observing that the $\rho$-action by a 3-braid leaves the third
coordinate invariant, i.e., $(f,g,h)*W(\sigma_1,\sigma_2) =
(f_1, g_1, h)$,
we have $(0,0,1)*\beta=(f,g, t^{4(k+l)})$ for some
$f,g \in \Z[t,t^{-1}]$.
Since $\rho(\beta)$ is trivial,
we obtain $k+l = 0$, which in turn implies that $\beta$ is
in $\langle \sigma_1, \sigma_2 \rangle = B_3  \subset B_4$.
The faithfulness of $\rho_3\otimes \Z_p$ by Lemma~\ref{thm:3braid}
leads to a contradiction.

\medskip
Now assume that $C$ contains the first two punctures
as Figure~\ref{fig:red}~(b).
The 4-braids represented by homeomorphisms that fix $C$
form a subgroup of $B_4$ generated by $\sigma_1$,
$x=\sigma_2\sigma_1^2\sigma_2$ and $y=\sigma_3$.
Since $\sigma_1$ commutes with both $x$ and $y$,
we write
$$\beta=\sigma_1^{k}W(x,y)$$
for an integer $k$ and a word $W$ on two letters.

By using the
relations $xyxy=yxyx$, $(xyxy)\sigma_1^2 = \Delta_4^2$ and
that $xyxy$ commutes with $x$, $y$ and $\sigma_1$;
we rewrite $\beta$ into another form by which we will
track $(0,0,1)*\beta$.

By substituting $x^{-1}$ by $(yxy)(xyxy)^{-1}$ and $y^{-1}$ by
$(xyx)(xyxy)^{-1}$ and then collecting $(xyxy)^{\pm1}$
to the left,
we have $W(x,y)=(xyxy)^{m}P(x,y)$ for some integer
$m$ and a positive word $P$ on two letters. We can assume that
we have moved $(xyxy)$ to the left as many as possible so that
neither $xyxy$ nor $yxyx$ occurs in $P$ as a subword.
We have
$$\beta=\sigma_1^{k}(xyxy)^{m}P(x,y)=\Delta_4^{2m}\sigma_1^{k-2m}P(x,y).$$

\medskip
We claim that $P$ contains both $x$ and $y$ as a subword.
If $P$ does not contain $y$, i.e., $P=x^l$ for some $\l\ge 0$,
then $\beta=\Delta_4^{2m}\sigma_1^{k-2m}x^l=
\Delta_4^{2m}\sigma_1^{k-2m}(\sigma_2\sigma_1^2\sigma_2)^l$
fixes the curve in Figure~\ref{fig:red}~(a).
By the previous argument $\beta$ is trivial.
If $P$ does not contain $x$, i.e.,
$P=y^l$ for some $l\ge 0$, then $\beta=\Delta_4^{2m}\sigma_1^{k-2m}y^l$.
From $(0,0,1)*\beta=(0,0,(-t)^{4m+l})$ and 
$(1,0,0)*\beta=((-t)^{4m+(k-2m)},0,0)$,
we have $l=-4m$ and $k=-2m$.
The exponent sum $e(\beta)=12m+(k-2m)+l=4m$ should equal zero
since $\rho(\beta)$ is trivial.
Therefore we have $m=l=k=0$, which implies that $\beta$ is trivial.

\medskip
Since $x$ and $y$ both
commutes with $\sigma_1$ and $\Delta_4^2$, by applying a conjugation
we may assume that
$P$ starts with $y$ and ends with $x$.
In Figure~\ref{fig:auto}~(a), we construct an automaton
that accepts a positive word  in $x,y$ without any occurrence of
$xyxy$ and $yxyx$.
Arbitrary paths  following the arrows give
words accepted by the automaton.
Now we have
$$\beta=\Delta_4^{2m}\sigma_1^{k-2m}Q(x,y,xy,yx,yxy,xyx)$$ for
some positive word $Q$ accepted by the automaton in
Figure~\ref{fig:auto}~(a).
Note that $Q$ starts by one of $\{y, yxy, yx\}$ and ends by
one of $\{x, xyx, yx\}$. In other words
$Q$ is represented by a path starting at the state
$Y$ and ending at the state $X$.

\medskip
We replace $xyx$ by $y^{-1}(xyxy)$,
$yxy$ by $x^{-1}(xyxy)$ and then collect all $(xyxy)$'s to the left
to obtain
$$\beta=\Delta_4^{2m_1}\sigma_1^{k_1}Q(x,y,xy,yx,x^{-1},y^{-1})$$
for some $k_1$ and $m_1$.

\begin{figure}
$$
{\includegraphics[scale=.9]{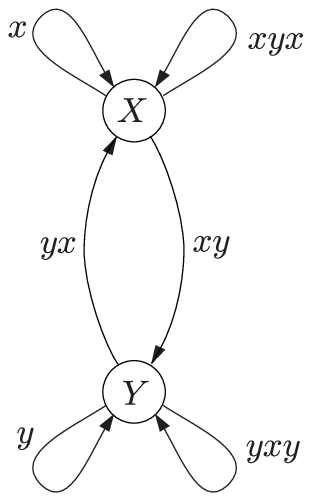}\atop\mbox{(a)}}\qquad\qquad
{\includegraphics[scale=.9]{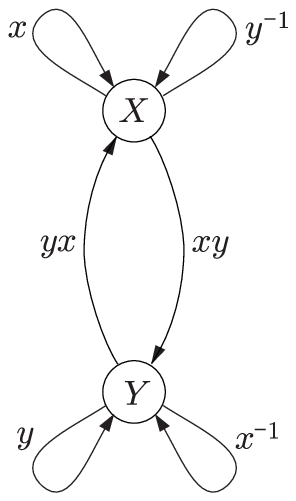}\atop\mbox{(b)}}
$$
\caption{Automata: (a) accepts all the words
that do not contain $xyxy$ or $yxyx$.}\label{fig:auto}
\end{figure}

\medskip
Consider the following two subsets of $\Z_p[t,t^{-1}]^3$.
\begin{eqnarray*}
V_X &=&\{(f,g,h)\in\Z_p[t,t^{-1}]^3\mid \deg g>\deg f,\ \deg g\ge\deg h\}\\
V_Y &=&\{(f,g,h)\in\Z_p[t,t^{-1}]^3\mid \deg h>\deg f,\ \deg h>\deg g\}
\end{eqnarray*}

The $\rho$-action of each arrow of the automaton
in Figure~\ref{fig:auto}~(b) is given as follows.
Let $\mathbf v=(f,g,h)\in\Z_p[t,t^{-1}]^3 $ be an arbitrary vector.
\begin{eqnarray*}
\mathbf v*x &=& (tf+(t^2-t)g+(1-t)h,\ t^3g+(1-t^2)h,\ h)\\
\mathbf v*y &=& (f,\ g,\ tg-th)\\
\mathbf v*(xy)&=& (tf+(t^2-t)g+(1-t)h,\ t^3g+(1-t^2)h,\ t^4g-t^3h)\\
\mathbf v*(yx)&=& (tf+(t^2-t)h,\ tg+(t^3-t)h,\ tg-th)\\
\mathbf v*x^{-1} &=&\textstyle
({t^{-1}}f+({t^{-3}}-{t^{-2}})g+({t^{-2}}-{t^{-3}})h,\
{t^{-3}}g+({t^{-2}}-{t^{-3}})h,\ h)\\
\mathbf v*y^{-1} &=&\textstyle (f,\ g,\ g-{t^{-1}}h)
\end{eqnarray*}
Then it is routine to check from the above formulae that
\begin{enumerate}
\item
$ V_X*x\subset V_X$,  $V_X*y^{-1}\subset V_X$
and $V_X*(xy)\subset V_Y$;
\item  $V_Y *y\subset V_Y$,  $V_Y*x^{-1}\subset
V_Y$ and $V_Y*(yx)\subset V_X$.
\end{enumerate}
Note that these are compatible with the automaton in
Figure~\ref{fig:auto}~(b).
If a path starts at $Y$ and ends at $X$
then the $\rho$-action of its braid word maps
$V_Y$ into $V_X$.
So we have $V_Y * Q \subset V_X$ for $Q=Q(x,y,xy,yx,x^{-1},y^{-1})$.

Since $(0,0,t^{4m_1}) \in V_Y$, we have
\begin{alignat*}{1}
(0,0,1)*\beta
&= (0,0,1)* \Delta_4^{2m_1}\sigma_1^{k_1} Q \\
&= (0,0,t^{4 m_1}) * \sigma_1^{k_1}  Q  \\
&= (0,0,t^{4 m_1}) * Q \quad \in V_X.
\end{alignat*}
Since $(0,0,1) \in V_Y$ and $V_X$ and $V_Y$ are disjoint,
$(0,0,1)*\beta \in V_X$ implies that
$\rho(\beta)$ is non-trivial.
\end{proof}

We remark that
the group generated by $x$ and $y$ is the
Artin group of Coxeter type $B_2$ and that $xyxy=yxyx$ is
the defining relation of the subgroup generated by $x$ and $y$.
So the subgroup generated by $x$, $y$ and $\sigma_1$ is
the direct product of the infinite cyclic subgroup
generated by $\sigma_1$ and the subgroup generated
by $x$ and $y$.

\let\s\sigma
\section{Non-Brunnian elements in $\ker \operatorname{Burau}(4)\otimes \Z_p$}
In~\cite{MR97m:20050} Cooper and Long obtain
a presentation of the image of $\rho_4\otimes \Z_2$.
As a corollary, $\rho_4\otimes \Z_2$ is not faithful.
In~\cite{MR99k:20077}, they compute
a presentation of a group containing the image of $\rho_4\otimes \Z_3$
as a finite index subgroup and
also give a non-trivial braid in the kernel explicitly.
In this section
we exhibit that the examples of Cooper and Long are not Brunnian.

\begin{figure}
\includegraphics[scale=.7]{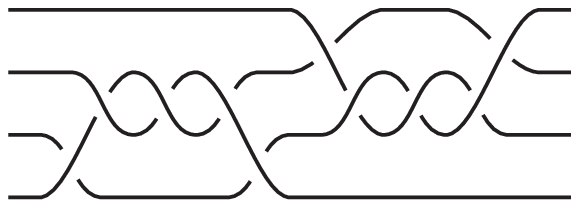}
\caption{The braid $\s_1^{-1}\s_2^3\s_1\s_3\s_2^{-3}\s_3^{-1}$,
whose fourth power is in the kernel of
$\operatorname{Burau}(4)\otimes \Z_2$.}
\label{fig:ex1}
\end{figure}

Let $\alpha_k=(\s_1^{-1}\s_2^k\s_1\s_3\s_2^{-k}\s_3^{-1})^4$ for $k\ne 0$.
See Figure~\ref{fig:ex1} for $\s_1^{-1}\s_2^3\s_1\s_3\s_2^{-3}\s_3^{-1}$.
The braid $\alpha_k$ comes from the fourth relation of Theorem 1.4
in~\cite{MR97m:20050} and
is in the kernel of $\beta_4\otimes \Z_2$.
$\alpha_k$ is not Brunnian because we obtain $\s_1^{4k}$
by forgetting the second and the fourth strands.

\begin{figure}
\includegraphics[scale=.7]{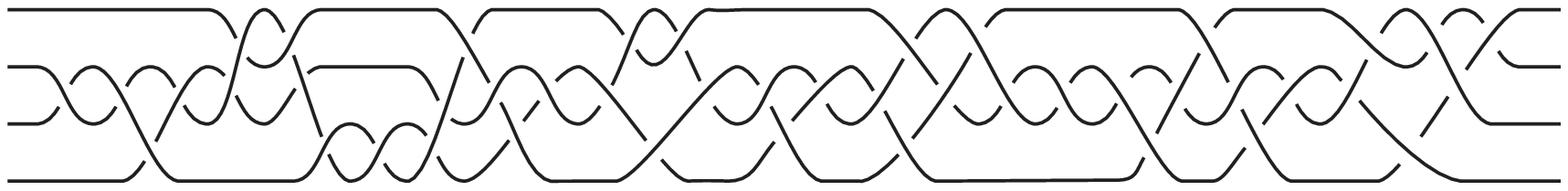}
\caption{A braid in the kernel of
$\operatorname{Burau}(4)\otimes \Z_3$}\label{fig:ex2}
\end{figure}

\begin{figure}
\includegraphics[scale=.7]{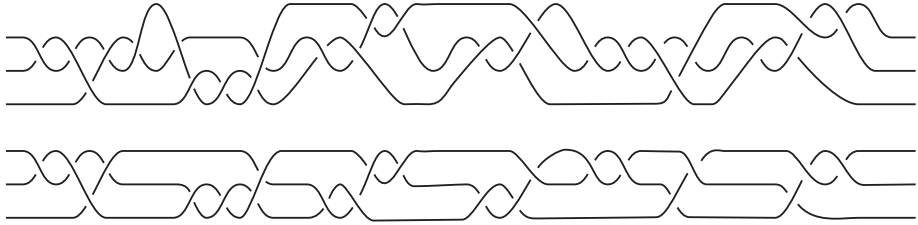}
\caption{Forgetting the fourth strand}\label{fig:ex3}
\end{figure}

Now let $\alpha$ be the braid
$$
\s_2^2\s_1\s_2^{-2}\s_3^{-2}
\s_2\s_1^{-3}\s_2^{-1}
\s_3\s_2^{-1}\s_1\s_2^2\s_3^{-2}\s_1^{-1}\s_2^{-2}
\s_1\s_2^{-2}\s_1\s_3\s_2^{-1}\s_3\s_2^3
\s_1\s_2^{-1}\s_3\s_2^{-1}\s_1\s_2^{-2}\s_1
\s_3^2\s_2\s_3^{-1}
$$
as Figure~\ref{fig:ex2}.
It is conjugate
to the braid given by \cite{MR99k:20077} as a non-trivial element of
$\ker \operatorname{Burau}(4)\otimes \Z_3$.
It is easy to see that $\alpha$ is not Brunnian.
If we forget the fourth strand from $\alpha$ as Figure~\ref{fig:ex3},
we get a non-trivial 3-braid
$$\alpha'=
\s_2^2\s_1\s_2^{-1}\s_1^{-3}\s_2^{-1}\s_1^2\s_2^{-2}
\s_1^{-2}\s_2^3\s_1^{-1}\s_2\s_1^{-1}\s_2^2
=(\s_2\s_1^{-1}\s_2\s_1^{-1}\s_2^2)^3\Delta_3^{-2}.
$$

\end{document}